\documentclass[10pt]{article}
\usepackage{pb-diagram}
\usepackage{amsmath, amsfonts}
\usepackage{amssymb}
\usepackage{amscd}

\newtheorem{Df}{Definition}[section]
\newtheorem{Te}[Df]{Theorem}
\newtheorem{Po}[Df]{Proposition}
\newtheorem{Cr}[Df]{Corollary}
\newtheorem{Lm}[Df]{Lemma}
\newtheorem{Ca}[Df]{Claim}
\newtheorem{Cn}[Df]{Conjecture}
\newtheorem{Ex}[Df]{Example}
\newtheorem{Rm}[Df]{Remark}

\newcommand{\Bdf}{\begin{Df}}
\newcommand{\Edf}{\end{Df}}
\newcommand{\Bte}{\begin{Te}}
\newcommand{\Ete}{\end{Te}}
\newcommand{\Bpo}{\begin{Po}}
\newcommand{\Epo}{\end{Po}}
\newcommand{\Bcr}{\begin{Cr}}
\newcommand{\Ecr}{\end{Cr}}
\newcommand{\Blm}{\begin{Lm}}
\newcommand{\Elm}{\end{Lm}}
\newcommand{\Bca}{\begin{Ca}}
\newcommand{\Eca}{\end{Ca}}
\newcommand{\Bcn}{\begin{Cn}}
\newcommand{\Ecn}{\end{Cn}}
\newcommand{\Bex}{\begin{Ex}}
\newcommand{\Eex}{\end{Ex}}
\newcommand{\Brm}{\begin{Rm}}
\newcommand{\Erm}{\end{Rm}}
\newcommand{\Bdm}{{\it Proof.}\ }
\newcommand{\Edm}{\rule{2mm}{2mm}}

\begin{document}

\title{\bf{A Koszul sign map}}
\author{Roland Berger}
\date{}

\maketitle

\begin{abstract}
We define a Koszul sign map encoding the Koszul sign convention. We view the Koszul sign map as a 2-cochain with respect to the cohomology of the permutation group. 
\end{abstract} 

\section{Introduction}

The Koszul sign convention is a rule which plays a fundamental role in graded algebra. This convention was used by Koszul in his thesis~\cite{ko:hom}. Actually, whenever the algebraic framework is based on graded categories as in homology theory~\cite{loday:cychom}, homotopy theory~\cite{fht:rht}, algebraic operads theory~\cite{lv:alop}, this convention is systematically used in the definitions of the new concepts naturally appearing in the theory, and in the statements of the new results as well. In general -- in particular in the cited books -- the convention is just stated, without further explanation. The aim of this paper is to propose a tool encoding precisely the Koszul sign convention.

In 1966, without any reference to Koszul, Boardman introduced a Principle of Signs~\cite{board:signs}, in order to make precise the way of inserting the right signs in the identities frequently obtained in algebraic topology. His approach lies on a class of various $n$-ary operations and on an involution, subjected to some axioms, allowing him to characterize the identities written in a standard form. Then he proved that the set of these identities is stable under natural algebraic transformations. Our approach is different, and it can be viewed as a preliminary step of Boardman's one. Although we fix a set of graded elements as in Boardman, we do not need to define the operations and the identities acting on these elements. For us, only the order of the elements inside of the result of an operation is important to produce a sign in front of the algebraic expression.

\section{A Koszul sign map}

The Koszul sign convention is used in various graded contexts. The objects on which the convention is applied are homogeneous, and the nature of the objects -- graded elements, graded maps -- depends on the context. However the convention does not depend on the nature of the objects, but just on their degrees. In our general setting, the homogeneous objects will be called \emph{symbols}. At each symbol, we associate a degree in $\mathbb{Z}$. 

Roughly speaking, the Koszul sign convention is the following: if in a manipulation of a monomial algebraic expression concerning on symbols naturally written from the left to the right, a symbol $f_i$ jumps over a symbol $f_j$ situated on the left or on the right of $f_i$, then the sign $(-1)^{|f_i||f_j|}$ appears in front of the expression, where $|f_i|$ denotes the degree of $f_i$. If the algebraic expression is a sum of monomial algebraic expressions, the convention is applied to each term.  

Let us note that the convention is independent of the algebraic operations included in the monomial algebraic expression. Although the manipulation passes from a monomial algebraic expression to another one, only the order of the objects in the initial and final expressions are significant.

For example, in the definition of a tensor product of graded linear maps
$$(f\otimes g)(a\otimes b)= (-1)^{|g||a|} f(a) \otimes g(b),$$
the ordered symbols in the initial -- final --  expression are
$f$, $g$, $a$, $b$ -- $f$, $a$, $g$, $b$.

In our setting, we are led to permute arbitrarly the symbols from the order of the initial expression. The result of a permutation might be organized in a final monomial algebraic expression respecting the final order, but such an expression would be irrelevant for us. 

Throughout the paper, we fix an integer $n\geq 2$, a sequence $f=(f_1, \ldots ,f_n)$ of symbols, and a sequence of degrees 
$$|f|=(|f_1|, \ldots ,|f_n|) \in \mathbb{Z}^n$$
called the degree of $f$.

We denote by $\mathcal{E}_f$ the set of the permutations of $f$. If $S_n$ is the group of permutations of $\{1, \ldots , n\}$, we have
$$\mathcal{E}_f= \{g=(g_1, \ldots ,g_n) \, ; \, \exists \rho \in S_n,\ g_i=f_{\rho^{-1}(i)}\ for \ i=1, \ldots, n\},$$
where $g_1, \ldots ,g_n$ are still seen as symbols. The degree of $g$ is defined by
$$|g|=(|g_1|, \ldots ,|g_n|)$$
with $|g_i|=|f_{\rho^{-1}(i)}|$. The group $S_n$ acts on the left on $\mathcal{E}_f$ by defining
$$\sigma (g)=(g_{\sigma^{-1}(1)}, \ldots ,g_{\sigma^{-1}(n)}),$$
for $g \in \mathcal{E}_f$ and $\sigma \in S_n$.

We want to define now a \emph{Koszul sign map}
\begin{eqnarray*}
\kappa : & S_n \times \mathcal{E}_f & \rightarrow \ \{\pm 1\}\\
         & (\sigma, g) & \mapsto \ \kappa (\sigma ,g),
\end{eqnarray*}
which respects the Koszul sign convention. For the moment, we state this convention in a rather heuristic form as follows.

\emph{If in the permutation $\rho : (g_1, \ldots ,g_n) \rightarrow (f_1, \ldots ,f_n)$ correcting the permutated sequence $(g_1, \ldots ,g_n)$ into the initial sequence $(f_1, \ldots ,f_n)$, $g_i$ jumps over $g_j$, then the sign $(-1)^{|g_i||g_j|}$ appears in $\kappa$.}

We begin to remark that it is not clear how to define a map
$$\kappa(-, f): S_n \rightarrow \{\pm 1\}$$
respecting the convention when $n=3$ -- if $n=2$, $\kappa(-, f): S_2 \rightarrow \{\pm 1\}$ is well-defined and is a group morphism.

In fact, after acting the transpositions $(1,2)$ and $(2,3)$ on $f=(f_1,f_2,f_3)$, we obtain
$$\kappa ((1,2),f)=(-1)^{|f_1||f_2|},$$
$$\kappa ((2,3),f)=(-1)^{|f_2||f_3|}.$$
Since $(2,3)(1,2)f=(f_3,f_1,f_2)$ is corrected into $(f_1,f_2,f_3)$ by jumping $f_3$ over $f_1$ and $f_2$, we obtain
$$\kappa ((2,3)(1,2),f)=(-1)^{|f_3|(|f_1|+|f_2|)},$$
so that $\kappa ((2,3)(1,2),f)\neq \kappa ((2,3),f)\, k((1,2),f)$ for a certain choice of the degrees. Therefore it is not possible to define a group morphism $\kappa(-, f): S_n \rightarrow \{\pm 1\}$ in great generality. However, if we set $g=(1,2)f=(f_2,f_1,f_3)$, we have $$\kappa ((2,3),g)=(-1)^{|g_2||g_3|}= (-1)^{|f_1||f_3|},$$
Thus we are led to the right formula
$$\kappa ((2,3)(1,2),f)= \kappa ((2,3),(1,2)(f))\, \kappa ((1,2),f).$$

In order to define $\kappa$ from this formula and its generalizations, we want to be sure that $\kappa (\sigma ,g)$ does not depend on the way to correct $g$ into $f$ by using transpositions. Consequently, we first define $\kappa$ on the free group $F_{n-1}$ generated by the transpositions $s_i=(i,i+1)$ for $i=1, \ldots n-1$. Let us recall that the group $S_n$ is defined by these generators and the following relations in $F_{n-1}$
\begin{equation} \label{relations}
s_i^2=e, \ \ s_is_j=s_js_i \ \mathrm{if} \ |i-j|>1, \ \ s_is_{i+1}s_i=s_{i+1}s_is_{i+1},
\end{equation}
where $e$ denotes the unit of the group $F_{n-1}$.

For $g$ in $\mathcal{E}_f$ and $1\leq i \leq n-1$, we set
\begin{equation} \label{action}
s_i(g)= s_i^{-1}(g)= (g_1, \ldots g_{i-1}, g_{i+1}, g_i, g_{i+2}, \ldots ,g_n),
\end{equation}
which defines an action $x(g)$ of the elements $x$ of the group $F_{n-1}$ on the set $\mathcal{E}_f$. This action induces naturally the action of $S_n$ on $\mathcal{E}_f$ defined above. 

We define the map
\begin{eqnarray} \label{freekappa}
\kappa : & F_{n-1} \times \mathcal{E}_f & \rightarrow \ \{\pm 1\}  \\
         & (x, g) & \mapsto \ \kappa (x , g) \nonumber
\end{eqnarray}
as follows. For any $g$ in $\mathcal{E}_f$, we set $ \kappa (e, g)=1$, and for 
$1\leq i \leq n-1$,
\begin{equation} \label{transposition}
\kappa (s_i, g)=\kappa (s_i^{-1}, g)=(-1)^{|g_i||g_{i+1}|}.
\end{equation}
Moreover, for any $x$ in $F_{n-1}$ decomposed in a reduced form
$$x =t_{i_1} \ldots t_{i_m}$$
where $t_{i_j}=s_{i_j}$ or $t_{i_j}=s_{i_j}^{-1}$, and $i_1, \ldots , i_m$ are in $\{1, \ldots , n-1\}$, we set
\begin{equation} \label{deffreekappa}
\kappa (x , g)=\kappa (t_{i_1}, t_{i_2} \dots t_{i_m} (g))\, \kappa (t_{i_2}, t_{i_3} \dots t_{i_m} (g)) \ldots \kappa (t_{i_m}, g).
\end{equation}
Reduced form means that two consecutive factors $t_{i_j}$ and $t_{i_{j+1}}$ are never inverse to each other. A reduced form being unique, the map (\ref{freekappa}) is well-defined. 

\Blm \label{lemma}
For $g$ in $\mathcal{E}_f$ and $1\leq i \leq n-1$, we have
\begin{equation} \label{inverses}
\kappa (s_i, s_i^{-1}(g))\, \kappa (s_i^{-1}, g)=1= \kappa (s_i^{-1}, s_i(g))\, \kappa (s_i, g).
\end{equation}
\Elm
\Bdm
Since $\kappa (s_i, g)=\kappa (s_i^{-1}, g)$ and $s_i(g)=s_i^{-1}(g)$, it suffices to verify the first equality. From (\ref{action}), we draw $\kappa (s_i, s_i^{-1}(g))=(-1)^{|g_{i+1}||g_i|}= \kappa (s_i^{-1}, g)$. \Edm 
\\

The lemma shows that the formula (\ref{deffreekappa}) extends to any decomposition $x =t_{i_1} \ldots t_{i_m}$ reduced or not. Therefore, one has
\begin{equation} \label{formula1}
\kappa (x \,y , g)=\kappa (x, y (g))\, \kappa (y, g)
\end{equation}
for $g$ in $\mathcal{E}_f$, $x$, $y$ in $F_{n-1}$, and consequently
\begin{equation} \label{formula2}
\kappa (x^{-1} , g)=\kappa (x, x^{-1} (g)).
\end{equation}

\Bpo \label{defkappa}
Passing through the relations (\ref{relations}), the map $\kappa : F_{n-1} \times \mathcal{E}_f \rightarrow \{\pm 1\}$ induces a map $\kappa : S_n \times \mathcal{E}_f \rightarrow \{\pm 1\}$.
\Epo
\Bdm
We will prove that, for each relation in (\ref{relations}) and for each fixed $g$, $\kappa (-,g)$ gives the same result on the left-hand side and on the right-hand side of the relation.

Fisrstly, $\kappa (s_i^2 , g)=\kappa (s_i, s_i(g))\, \kappa (s_i, g)$ according to (\ref{formula1}), hence $\kappa (s_i^2 , g)=1$ by the lemma.

Secondly, let us suppose that $j>i+1$, so that
$$s_j(g)= (g_1, \ldots g_{i}, g_{i+1}, \ldots g_{j+1}, g_{j}, \ldots ,g_n),$$
thus $\kappa (s_i, s_j(g))=(-1)^{|g_i||g_{i+1}|}$ and $\kappa (s_j, g)=(-1)^{|g_j||g_{j+1}|}$. Using (\ref{formula1}), we obtain
$$\kappa (s_i s_j, g)=(-1)^{|g_i||g_{i+1}|+|g_j||g_{j+1}|}.$$
The same if $j<i-1$. So, if $|i-j|>1$, we obtain the expected equality
$$\kappa (s_i s_j, g)= \kappa (s_j s_i, g).$$

Thirdly, from (\ref{formula1}), we draw
\begin{equation} \label{braid}
\kappa (s_i \,s_{i+1} \, s_i , g)=\kappa (s_i, s_{i+1} \, s_i (g))\, \kappa (s_{i+1}, s_i (g))\,\kappa (s_i, g).
\end{equation}
Using $s_i(g)= (g_1, \ldots g_{i+1}, g_i, g_{i+2}, \ldots ,g_n)$ and $s_{i+1} s_i(g)= (g_1, \ldots g_{i+1}, g_{i+2}, g_i, \ldots ,g_n)$, the formula (\ref{braid}) implies
$$\kappa (s_i \,s_{i+1} \, s_i , g)=(-1)^{|g_{i+1}||g_{i+2}|+|g_i||g_{i+2}|+|g_i||g_{i+1}|}.$$
Using $s_{i+1}(g)= (g_1, \ldots g_i, g_{i+2}, g_{i+1}, \ldots ,g_n)$ and $s_i s_{i+1}(g)= (g_1, \ldots g_{i+2}, g_i, g_{i+1}, \ldots ,g_n)$, the formula (\ref{braid}) in which $i$ and $i+1$ are exchanged gives
$$\kappa (s_{i+1} \, s_i \,s_{i+1} , g)=(-1)^{|g_i||g_{i+1}|+ |g_i||g_{i+2}| + |g_{i+1}||g_{i+2}|}.$$
Then we arrive to  $\kappa (s_i \,s_{i+1} \, s_i , g)=\kappa (s_{i+1} \, s_i \,s_{i+1} , g)$. 

An equivalent way to say what we have obtained is the following. Writing each relation (\ref{relations}) as an equality $r=e$ for a certain element $r$ in $F_{n-1}$, we have $\kappa (r,g)=1$ for any $g$ -- it suffices to apply the formula (\ref{formula1}).

Now for any $x$ in $F_{n-1}$, we have
\begin{eqnarray*}
\kappa (x r x^{-1},g) & = & \kappa (x, r (x^{-1} (g))\, \kappa (r, x^{-1}(g))\, \kappa (x^{-1},g) \\
                           & = & \kappa (x, x^{-1} (g))\, \kappa (x^{-1},g)
\end{eqnarray*}
since $r(g')=g'$ for any $g'$ in $\mathcal{E}_f$. Therefore (\ref{formula2}) implies that $\kappa (x r x^{-1},g)=1$ for any $g$ in $\mathcal{E}_f$.

Using again (\ref{formula1}), we deduce $\kappa( x r x^{-1} x' r' x'^{-1},g)=1$ for any $g$ in $\mathcal{E}_f$, and any elements $r$, $r'$ of $F_{n-1}$ generating the relations. Inductively, we obtain that $\kappa (a,g)=1$ for any $g$ in $\mathcal{E}_f$ and any $a$ in the normal subgroup of $F_{n-1}$ generated by the elements $r$. More generally, $\kappa (a x,g)=\kappa (x,g)$ for any $a$ as previously and any $x$ in $F_{n-1}$.
\Edm
\\

Our construction of the map $\kappa : S_n \times \mathcal{E}_f \rightarrow \{\pm 1\}$ shows immediately the following proposition. This proposition could be used as a definition. 
\Bpo \label{caractkappa}
The map $\kappa : S_n \times \mathcal{E}_f \rightarrow \{\pm 1\}$ is the unique map such that
\begin{enumerate}
\item $\forall \sigma \in S_n,\, \forall \tau \in S_n,\, \forall g \in \mathcal{E}_f, \ \kappa (\sigma \,\tau , g)=\kappa (\sigma, \tau (g))\, \kappa (\tau, g)$,
\item $\forall i \in \{1, \ldots , n-1\},\, \forall g \in \mathcal{E}_f, \ \kappa (s_i, g)=(-1)^{|g_i||g_{i+1}|}$.
\end{enumerate}
\Epo

From $f$ and its degree $|f|$, we have constructed the map $\kappa$ which should be rather denoted by $\kappa_f$.
For any $f'$ in $\mathcal{E}_f$, one has $\mathcal{E}_{f'}=\mathcal{E}_f$. The degree of $f'$ is obviously defined from the degree of $f$. Then Proposition \ref{caractkappa} shows that the map
$$\kappa_{f'} : S_n \times \mathcal{E}_{f'} \rightarrow \{\pm 1\}$$
coincides with $\kappa _f$. However, it is possible that
$$\kappa _f(-,f) \neq \kappa _{f'}(-,f'),$$
as in the example $f=(f_1,f_2,f_3)$ and $f'=(f_1,f_3,f_2)$, with $|f_1|$ and $|f_2|$ odd, $|f_3|$ even.

\Bex
\emph{Take $n=5$ and $g=(f_4,f_1,f_3,f_5,f_2)$, so that $g=\rho (f)$ where}
$$\rho=\left( \begin{array}{ccccc} 4 & 1 & 3 & 5 & 2 \\
  1 & 2 & 3 & 4 & 5 \end{array} \right).$$
\emph{Using a bubble sort, we find that $\kappa(\rho ,f)=(-1)^Z$ with $z_i=|f_i|$ and}
$$Z=z_1z_4 + z_3z_4 + z_2z_5 + z_2z_4 + z_2z_3.$$
\Eex

\setcounter{equation}{0}

\section{When $\kappa (-,g)$ is a group morphism}

An integer $n\geq 2$, a sequence $f=(f_1, \ldots ,f_n)$ of symbols, and a sequence of degrees $|f|=(|f_1|, \ldots ,|f_n|) \in \mathbb{Z}^n$ being given, we have defined the map $\kappa =\kappa_f$ in the previous section. If $n=2$, then $\kappa (-, (f_1,f_2))=\kappa (-, (f_2,f_1))$ is always a group morphism from $S_2$ to $\{\pm 1\}$.

Let us suppose that $n\geq 3$. We want to know when
$$\kappa (-,g) : S_n  \rightarrow \{\pm 1\}$$
is a group morphism. From 1. in Proposition \ref{caractkappa}, it is the case if and only if for any $\sigma$ and $\tau$ in $S_n$, one has
$$\kappa (\sigma, \tau (g))= \kappa (\sigma, g),$$
that is, if and only if, for any $\sigma$ in $S_n$ and $h$ in $\mathcal{E}_f$, 
$$\kappa (\sigma, h)= \kappa (\sigma, g),$$
which implies that $\kappa (- , h)$ is a group morphism as well.

So, it suffices to examine when $\kappa (- , f)$ is a group morphism, and we have seen that it is the case if and only if
\begin{equation} \label{morphism1}
\forall \sigma \in S_n,\, \forall \tau \in S_n,\, \forall g \in \mathcal{E}_f, \ \kappa (\sigma, \tau (g))=\kappa (\sigma, g).
\end{equation}
From 1. in Proposition \ref{caractkappa}, it is equivalent to 
\begin{equation} \label{morphism2}
\forall i \in \{1, \ldots , n-1\}, \, \forall j \in \{1, \ldots , n-1\}, \, \forall g \in \mathcal{E}_f, \ \kappa (s_i, s_j(g))=\kappa (s_i, g).
\end{equation}

It is clear that (\ref{morphism2}) holds whenever $i=j$ or $|i-j|>1$. Thus (\ref{morphism2}) is equivalent to
\begin{equation} \label{morphism3}
\forall i \in \{1, \ldots , n-2\}, \, \forall g \in \mathcal{E}_f, \ \kappa (s_i, s_{i+1}(g))=\kappa (s_i, g) \ \mathrm{and} \ \kappa (s_{i+1}, s_i(g))=\kappa (s_{i+1}, g).
\end{equation}
Some calculations included in the proof of Proposition \ref{defkappa} show that it is equivalent to
\begin{equation} \label{morphism4}
\forall i \in \{1, \ldots , n-2\}, \, \forall g \in \mathcal{E}_f, \ (-1)^{|g_i||g_{i+2}|} = (-1)^{ |g_i||g_{i+1}|} = (-1)^{|g_{i+1}||g_{i+2}|},
\end{equation}
which, in turn, is equivalent to say that among the parities of the triplets $(|g_i|,|g_{i+1}|,|g_{i+2}|)$, only the cases of one even parity and two odd parities are forbidden.

The last condition is satisfied if all the degrees $|f_1|, \ldots ,|f_n|$ have the same parities, otherwise if only one is odd. Conversely, if the last condition is satisfied and if $|f_1|, \ldots ,|f_n|$ do not have the same parities, it is not possible to find one even parity and two odd parities among them -- a suitable permutation of $f$ then providing a forbidden triplet. We have obtained the following.

\Bpo \label{kappamorphism}
Suppose that $n\geq 2$, $f=(f_1, \ldots ,f_n)$, and $|f|=(|f_1|, \ldots ,|f_n|) \in \mathbb{Z}^n$ are given. 
\begin{enumerate}
\item The map $\kappa (-,g) : S_n \rightarrow \{\pm 1\}$ is a group morphism for an element $g$ of $\mathcal{E}_f$ if and only if $\kappa (-,f)$ is a group morphism, and in this case, all the group morphisms $\kappa (-,g)$ are equal. 
\item The map $\kappa (-,f)$ is a group morphism if and only if either all the integers $|f_1|, \ldots ,|f_n|$ have the same parities or only one is odd among them.
\item The map $\kappa (-,f)$ is constant equal to $1$ if and only if all the integers $|f_1|, \ldots ,|f_n|$ are even, with possibly one exception.  
\end{enumerate}
\Epo

\setcounter{equation}{0}

\section{A cohomological interpretation of $\kappa_f$}

An integer $n\geq 2$, a sequence $f=(f_1, \ldots ,f_n)$ of symbols, and a sequence of degrees $|f|=(|f_1|, \ldots ,|f_n|) \in \mathbb{Z}^n$ are given. For $\sigma$ and $\rho$ in $S_n$, we put
$$c_f(\sigma, \rho) = \kappa_f (\sigma, \rho(f)).$$
So we define a map $c_f : S_n \times S_n \rightarrow \{\pm 1\}$. This is the unique map such that
\begin{enumerate}
\item $\forall \sigma \in S_n,\, \forall \tau \in S_n,\, \forall \rho \in S_n, \ c_f (\sigma \,\tau , \rho)=c_f (\sigma, \tau \, \rho)\, c_f (\tau, \rho)$,
\item $\forall i \in \{1, \ldots , n-1\},\, \forall \rho \in S_n, \ c_f (s_i, \rho)=(-1)^{|f_{\rho^{-1}(i)}||f_{\rho^{-1}(i+1)}|}$.
\end{enumerate}

We want to regard the map $c_f$ as a 2-cochain for the cohomology of the group $S_n$ with coefficients in the multiplicative group $\{\pm 1\}$. The automorphism group of $\{\pm 1\}$ is formed of $+Id$ and $-Id$, and it is identified to the group $\{\pm 1\}$. Then a structure of $S_n$-module on the group $\{\pm 1\}$ is equivalent to the datum of a group morphism
$$u : S_n \rightarrow \{\pm 1\},$$
the action of $\sigma \in S_n$ on $+1$ or $-1$ being given by the product by $u(\sigma)$. It is well-known that there are only two such morphisms $u$: the constant morphim $u=1$ and the signature $u=sgn$. 

Let us choose a structure $u$ of $S_n$-module on $\{\pm 1\}$. Let us calculate the coboundary operator $\delta$ of the 2-cochain $c_f$. For $\sigma$, $\tau$ and $\rho$ in $S_n$, one has
$$\delta (c_f)(\sigma, \tau, \rho)= (\sigma . c_f(\tau, \rho))\, c_f(\sigma \tau, \rho)^{-1}\, c_f(\sigma, \tau \rho) \, c_f(\sigma, \tau)^{-1}.$$
Since $\sigma . c_f(\tau, \rho)= u(\sigma)\, c_f(\tau, \rho)$, the relation 1. just above shows that
$$\delta (c_f)(\sigma, \tau, \rho)=  u(\sigma)\, c_f(\sigma, \tau)^{-1}.$$
Thus $c_f$ is a 2-cocycle if and only if $c_f(\sigma, \tau)= u(\sigma)$ for any $\sigma$ and $\tau$. This condition implies that $c_f (-, \tau): S_n \rightarrow \{\pm 1\}$ is a group morphism for any $\tau$, thus either all the integers $|f_1|, \ldots ,|f_n|$ have the same parities or only one is odd among them (Proposition \ref{kappamorphism}).

Conversely, if either all the integers $|f_1|, \ldots ,|f_n|$ have the same parities or only one is odd among them, then the numbers $u_i=(-1)^{|f_{i}||f_{i+1}|}$ are all equal. By $u(s_i)= u_i$ for $i = 1, \ldots , n-1$, we define a group morphism $u : S_n \rightarrow \{\pm 1\}$, hence a structure of $S_n$-module on $\{\pm 1\}$ for which $c_f$ is a 2-cocycle. Moreover we have $u=c_f(-, e)$, so that
$$\delta (u)(\sigma, \tau)= u(\sigma) c_f(\tau, e)\, c_f(\sigma \tau, e)\, c_f(\sigma, e).$$
Then the relation 1. just above implies that $\delta (u)=c_f$. Let us sum up what we have obtained.

\Bpo \label{2-cocycle}
Let us suppose that $n\geq 2$, $f=(f_1, \ldots ,f_n)$, and $|f|=(|f_1|, \ldots ,|f_n|) \in \mathbb{Z}^n$ are given. Let us endow the group $\{\pm 1\}$ with the $S_n$-module structure defined by a group morphism $u : S_n \rightarrow \{\pm 1\}$. Then the map
$$c_f : S_n \times S_n \rightarrow \{\pm 1\}$$
is a 2-cocycle with coefficients in the group $\{\pm 1\}$ if and only if either all the integers $|f_1|, \ldots ,|f_n|$ have the same parities or only one is odd among them.

When this condition holds, $u=1$ if and only if all the integers $|f_1|, \ldots ,|f_n|$ are even, with possibly one exception, and $u=sgn$ if and only if all the integers $|f_1|, \ldots ,|f_n|$ are odd. Moreover, in both cases, $c_f$ is equal to the coboundary of $u$. 
\Epo

The 2-cochain $c_f$ is not symmetric in general. In fact, $c_f(e,\tau)=1$ for any $\tau$, while $c_f(s_i, e)= (-1)^{|f_{i}||f_{i+1}|}$.

\textbf{Question.} Find another cohomological interpretation for which $c_f$ is always a 2-cocycle, and is a 2-coboundary if and only if $c_f(-,e)$ is a group morphism.

\vspace{0.5 cm} \textsf{Roland Berger: Univ Lyon, UJM-Saint-\'Etienne, CNRS UMR 5208, Institut Camille Jordan, F-42023, Saint-\'Etienne, France}

\emph{roland.berger@univ-st-etienne.fr}\\

\end{document}